\documentclass{amsart}

\usepackage{amsthm}
\usepackage{graphicx}
\usepackage{amssymb}
\usepackage{psfrag}

\newtheorem{teo}{Theorem}[section]
\newtheorem{lemma}[teo]{Lemma}
\newtheorem{prop}[teo]{Proposition}
\newtheorem{cor}[teo]{Corollary}
\newtheorem{conj}[teo]{Conjecture}
\newtheorem{question}[teo]{Question}

\theoremstyle{definition}
\newtheorem{defi}[teo]{Definition}
\newtheorem{rem}[teo]{Remark}
\newtheorem{example}[teo]{Example}

\theoremstyle{remark}
\newtheorem{prof}[teo]{Proof of}

\def\mn{\mathbb{N}}
\def\mc{\mathbb{C}}
\def\mr{\mathbb{R}}
\def\mz{\mathbb{Z}}
\def\ns{\negmedspace}
\def\nns{\negthickspace}

\begin{document}
\title{Colored Jones invariants of links in $S^3\#_kS^2\times S^1$ and the Volume Conjecture}

\author[Costantino]{Francesco Costantino}
\address{Institut de Recherche Math\'ematique Avanc\'ee\\
  Rue Ren\'e Descartes 7\\
  67084 Strasbourg, France}
\email{f.costantino@sns.it}


\begin{abstract}
We extend the definition of the colored Jones polynomials to framed links and trivalent graphs in $S^3\#_kS^2\times S^1$ using a state-sum formulation based on Turaev's shadows.
Then, we prove that the natural extension of the Volume Conjecture is true for an infinite family of hyperbolic links.
\end{abstract}

\maketitle

\tableofcontents
\section{Introduction}
Since its discovery in the early eighties, the Jones polynomial of links in $S^3$ has been one of the most studied objects in low-dimensional topology. Despite this, at present, it is not yet completely clear which topological informations are carried and how are encoded by this invariant and more in general by the larger family of Quantum Invariants.
One of the most important conjectural relations between the topology of a link in a manifold and its Quantum Invariants has been given by Kashaev through his Volume Conjecture for hyperbolic links $L$ in $S^3$ (\cite{K}), based on complex valued link invariants $<L>_d$ constructed by using planar $(1,1)$-tangle presentations of $L$ and constant Kashaev's $R$-matrices (\cite{K2}).
Later Murakami-Murakami (\cite{MM}) identified Kashaev's invariants as special evaluations of certain colored Jones polynomials, and extended the Volume Conjecture to non-necessarily hyperbolic links in $S^3$ by replacing the hyperbolic volume of $S^3\setminus L$ with the (normalized) Gromov norm.
Let us state more precisely this Volume Conjecture. Let us recall that, using the irreducible representations of $U_q(sl_2(\mc))$, it is possible to assign to each link $L\subset S^3$ a sequence of Laurent polynomials $J_d(L),\ d\geq 2$, called the $d$-colored Jones polynomials of $L$ so that $J_d(L)$ only depends on the isotopy class of $L$ and all the polynomials associated to the unknot are equal to $1$ (see \cite{J} and \cite{Tu3}). 
The Volume Conjecture states the following:
$$\lim_{d\to \infty} \frac{1}{d}log(|ev_d(J_d(L))|)=\frac{Vol(L)}{2\pi}$$
where by $Vol(L)$ we denote the (normalized) Gromov norm of $S^3-L$ and by $ev_d$ the evaluation in $e^{\frac{2\pi\sqrt{-1}}{d}}$.
The VC in this form has been formally checked for the Figure Eight knot (\cite{MM}) and for a class of non-hyperbolic knots including torus knots (\cite{KT}, \cite{MM3}). Moreover, there are experimental evidences of its validity for knots $6_3$, $8_9$ and $8_{20}$ and for the Whitehead link (\cite{MMOTY}).

More recently Baseilhac and Benedetti have constructed a $(2+1)$-dimensional so called ``quantum {\it hyperbolic} field theory" (QHFT), which includes invariants $H_d(W,L,\rho)$ ($d\geq 1$ being any odd integer), where $W$ is an arbitrary oriented closed $3$-manifold, and either $L$ is a non-empty link in $W$ and $\rho$ is a principal flat $PSL(2,\mc)$-bundle on $W$ (up to gauge transformation) (\cite{BB1}, \cite{BB2}), or $L$ is a {\it framed} link and $\rho$ is defined on $W\setminus L$ (\cite{BB3}). The state-sums giving $H_d(W,L,\rho)$ for $d>1$ are in fact non commutative generalizations (see \cite{BB2}) of known simplicial formulas (\cite{N}) for $CS(\rho)+i Vol(\rho)$, the Chern-Simons invariant and the volume of $\rho$, which coincide with the usual geometric ones when $\rho$ corresponds to the holonomy of a finite volume hyperbolic $3$-manifold. Each $H_d(W,L,\rho)$ is well defined up to $dth$-roots of unity multiplicative factors. It is a non trivial fact (see \cite{K2}) that when $W=S^3$ and $\rho$ is necessarily  the trivial flat $PSL(2,\mc)$-bundle on $S^3$, the Kashaev's invariants $<L>_d$ coincide with $H_d(S^3,L,\rho_0)$ up to the above phase ambiguity. We recall that in \cite{BB1}, \cite{BB2}, \cite{BB3} different instances of general ``Volume Conjectures" have been formulated in the QHFT framework. Roughly speaking, these predict that, in suitable geometric situations (for example related to hyperbolic Dehn filling and the convergence of closed hyperbolic $3$-manifolds to manifolds with cusps), the quantum state sums asymptotically recover a classical simplicial formula, when $d\rightarrow \infty$. However, we note that Kashaev's VC (even if reformulated in terms of $H_d(S^3,L,\rho_0)$) is {\it not} a specialization of such QHFT Volume Conjectures.

One of the many difficulties in checking the Volume Conjecture in any of its instances is that in general the formulas describing the $d$-colored Jones polynomials of a link in $S^3$, or the state-sums giving $H_d(W,L,\rho)$, become more and more complicated while $d$ grows. Moreover, both diagrams of hyperbolic knots in $S^3$ and decorated triangulations of $(W,L)$ supporting the QHFT state-sums are usually quite complicated, so that, in the hyperbolic case, one almost immediately faces ugly formulas.

In order to test the ``VC" in new hyperbolic cases, we have by-passed this problem in the following way:
\begin{enumerate}
\item There is a particularly interesting {\it infinite} family of hyperbolic links in $\#_k S^2\times S^1$, $k\geq 2$, called {\it universal hyperbolic links}, discovered by the author and D.P. Thurston (\cite{CT2});
\item each universal hyperbolic link has a particularly simple presentations through Turaev's shadows theory;
\item for links in $S^3$, there is a state-sum formulation of $d$-colored Jones polynomials based on Turaev's shadows (\cite{KR},\cite{Tush});
\item then we extend the shadow definition of the $d$-colored Jones polynomials to links in $\#_k S^2\times S^1$, $k\geq 1$;
\item finally we prove that a natural extension of Kashaev's VC, in terms of these generalized Jones invariants, actually holds for universal hyperbolic links. 
\end{enumerate}

A natural question arises whose answer is still unknown to us:
\begin{question} Are the evaluations of the Jones invariants for a link $L$ in $\#_k S^2\times S^1$ (for odd $d$), entering this extended Kashaev's VC, special instances of the above QHFT invariants (presumably $H_d(\#_k S^2\times S^1,L,\rho_0)$, $\rho_0$ being again the trivial flat $PSL(2,\mc)$-bundle?
\end{question}
We now briefly summarize the contents of the paper. After recalling (in Section \ref{simplepoly}) the basic machinery and results on shadows of links in $3$-manifolds, in Section \ref{sec:jones}, we extend the definition of the colored Jones invariants using the shadow-based state-sum formulation of the Reshetikhin-Turaev invariants of $3$-manifolds (see \cite{Tush} and \cite{KR}). As a result we prove Theorem \ref{teo:jones} restated below in a slightly simplified form:
\begin{teo}
Let $L$ be a framed trivalent graph in an oriented $3$-manifold $N$ diffeomorphic to a connected sum of $k$ copies $S^2\times S^1$ or to $S^3$. There exist rational functions $J_d(N,L),\ d\geq 2$ which are invariant up to homeomorphism of the pair $(N,L)$ and extend the Reshetikhin-Turaev invariants. Moreover, if $N=S^3$ and $L$ is a link, $J_d(N,L)$ is the $d$-colored framed Jones polynomial normalized so that its value on the $0$ framed unknot is $1$. 
\end{teo}
For the sake of completeness, let us note that there exists also a generalization of Jones polynomials to the case of links in lens spaces provided by Hoste and Przytycki (\cite{HP}).

In the same section, we give examples of computations and recall the construction of the universal hyperbolic links. After recalling the many interesting properties shared by these links, we provide a particularly simple formula for their $d$-colored Jones invariants.

In the last section, after discussing how to extend the VC to the case of links in $S^3\#_k S^2\times S^1$ through the $d$-colored Jones invariants and recalling some definitions and results on the Lobatchevskji function, we prove the following (Theorem \ref{teo:vcvera} below):
\begin{teo}
The Extended Volume Conjecture is true for all the universal hyperbolic links.
\end{teo}
  
{\bf Acknowledgments.} 
I wish to express my gratitude to Stephane Baseilhac, Riccardo Benedetti, Dylan Thurston and Vladimir Turaev for their illuminating suggestions and comments. I also wish to warmly thank Manolo Eminenti and Nicola Gigli for our crucial ``home discussions".
\section{Recalls on shadows of graphs and state sums}\label{simplepoly}
In this section we recall the basic definitions and facts about shadows; no new result is proved. For a more detailed account, see \cite{Tu} and \cite{Coin}.
\subsection{Shadows of graphs}\label{sec1}
\begin{figure}
  \centerline{\includegraphics[width=8.4cm]{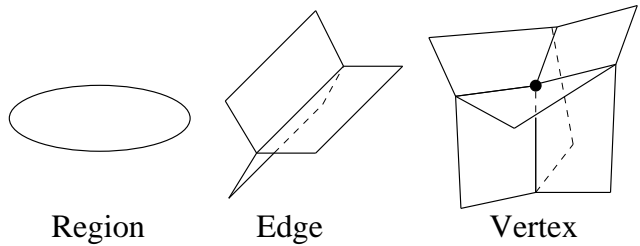}}
  \caption{The three local models of a simple polyhedron. }
  \label{fig:singularityinspine}
\end{figure}

\begin{figure}
  \centerline{\includegraphics[width=5.4cm]{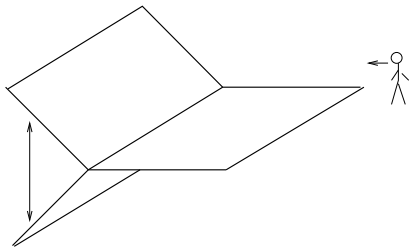}}
  \caption{The picture sketches the position of the polyhedron in a
    $3$-dimensional slice of the ambient $4$-manifold.  The direction
    indicated by the vertical double arrow is the one along which the
    two regions touching the horizontal one get separated.  }
  \label{fig:divergingdirection}
\end{figure}

\begin{defi} [Simple polyhedron] A simple polyhedron $P$ is a $2$-dimensional CW complex whose
local models are those depicted in Figure
\ref{fig:singularityinspine}; the set of points whose neighborhoods have models
of the two rightmost types is a $4$-valent graph, called
{\it singular set} of the polyhedron and denoted by $Sing(P)$. The
connected components of $P-Sing(P)$ are the {\it regions} of $P$.
The set of points of $P$ whose local models correspond to the boundaries of the blocks shown in the figure is called the {\it boundary} of $P$ and is denoted by $\partial P$; $P$ is said to be {\it closed} if $\partial P=\emptyset$. If a region contains no edges of $\partial P$ is called {\it internal}, otherwise  {\it external}.
The {\it complexity} $c(P)$ of a simple polyhedron $P$ is its number of vertices.
\end{defi}
From now on all the manifolds will be smooth and oriented unless explicitly stated and all the  polyhedra will be simple.
\begin{defi}
Let $M$ be a compact $4$-manifold with boundary and $P$ be a polyhedron embedded in $M$ so that $\partial P\subset \partial M$. We say that $P$ is a \emph{shadow} of the pair $(\partial M,\partial P)$ if the following holds:
\begin{enumerate}
\item $P$ is flat in $M$, i.e. $\forall p\in P$ there exists a local chart $U$ of $M$ centered in $p$ such that $(U,U\cap P)$ appears as in Figure \ref{fig:singularityinspine} and, in particular, $U\cap P$ is contained in a $3$-dimensional smooth ball in $M$.  
\item $M-P$ is diffeomorphic to $(\partial M-\partial P)\times [0,1)$ (i.e. $M$ collapses on $P$).
\end{enumerate} 
\end{defi}
If $P$ is a surface with framed boundary in $\partial M$, we can define an integer self-intersection number for $P$ in $M$. More in general, recalling that a framed trivalent graph in a $3$-manifold is the pair of a graph and a surface with boundary collapsing on it, the following holds:
\begin{prop}\label{prop:gleam}
Let $M$ be a compact $4$-manifold and $P$ be a flat polyhedron in $M$. If $\partial P\subset \partial M$ is framed in $\partial M$, then there exists a canonical coloring of the regions of $P$ with half-integers called \emph{gleam} induced by the embedding of $P$ in $M$ and by the framing on $\partial P$.
\end{prop}
\begin{prf}{1}{Let us give a sketch of the proof, for a detailed account see \cite{Tu}.
Let $D$ be a region of $P$, $\overline{D}$ be
the natural compactification of the (open) surface represented by
$D$ and $cl(D)$ the closure of $D$ in $P$. The embedding of $D$ in $P$ extends to a map
$i:\overline{D}\to cl(D)$ which is injective in $int(\overline{D})$,
locally injective on $\partial \overline{D}$ and which sends
$\partial \overline{D}$ into $Sing(P)\cup \partial P$. 

Orient $D$ arbitrarily and orient its normal bundle $n$ in $M$ so that the orientation on the global space of $n$ coincides with that of $M$.
Pulling back $n$ to $\overline{D}$ through $i$, we get an oriented disc bundle over $\overline{D}$ we will call the ``normal bundle" of $\overline{D}$; we claim that the projectivization of this bundle comes equipped with a section defined on $\partial \overline{D}$. 
Indeed, since $P$ is locally flat in $M$, for each point $p\in \partial \overline{D}$, if $i(p)\in Sing(P)$ there exists a smooth $3$-ball $B_{i(p)}\subset M$ around $i(p)$ in which $P$ appears as in Figure \ref{fig:divergingdirection}. Then, the intersection in $T_{i(p)} M$ of the ($3$-dimensional) tangent bundle of $B_{i(p)}$ with the normal bundle to $cl(D)$ in $i(p)$ gives a normal direction to $cl(D)$ in $i(p)$ (indicated in the figure) whose pull-back through $i$ is the seeked section in $p$; if $i(p)\in \partial P$ we use the framing on $\partial P$ to define the normal line to $cl(D)$ in $i(p)$ and then pull-it back through $i$. Hence a section of the projectivized normal bundle of $\overline{D}$ is defined on all $\partial \overline{D}$: we then define $gl(D)$ be equal to $\frac{1}{2}$ times the obstruction to extend this section to the whole $\overline{D}$; such an obstruction is an element of $H^2(\overline{D},\partial \overline{D};\pi_1(S^1))$, which is canonically identified with $\mathbb{Z}$ since $M$ is oriented. }\end{prf}
\begin{rem}\label{rem:framing}
We stress here that our definition of framing on a trivalent graph implies that, when the graph is a circle, then the framing can also be given by a Moebius band. 
\end{rem}
Proposition \ref{prop:gleam} shows that a flat embedding of a polyhedron in a $4$-manifold naturally equips the polyhedron with a combinatorial datum, encoded by the gleam. The following result is a converse: 
\begin{teo}[Turaev, \cite{Tu}]\label{teo:reconstruction}
Let $(P,gl)$ be a simple polyhedron equipped with gleams. It is possible to canonically thicken $(P,gl)$ to a smooth, compact and oriented $4$-manifold $M_{(P,gl)}$ such that
$P$ is a shadow of $M$ and $\partial P$ is framed in $\partial M$. 
Moreover, if $P$ is embedded in a $4$-manifold $M$, and $gl$ is the gleam induced on $P$ by its embedding (see Proposition \ref{prop:gleam}) then $M_{(P,gl)}$ is diffeomorphic to a neighborhood of $P$ in $M$.
\end{teo}\label{teo:thickening}
\begin{rem}
Our formulation of Theorem \ref{teo:reconstruction} is slightly imprecise: indeed a compatibility condition has to be satisfied on the gleams of the internal regions of $P$ in order for the thickening to be feasible (see \cite{Tu} and \cite{Coin}). However, this does not cause any problems to our following arguments since in our applications the polyhedra will contain no internal regions.
\end{rem}
\begin{defi}
Let $N$ be a closed $3$-manifold and $T$ be a framed trivalent graph in $N$ (possibly without vertices or even empty).
A pair $(P,gl)$ where $P$ is a polyhedron and $gl$ is a gleam on it, is a \emph{shadow} of $(N,T)$ if the pair $(\partial M_{(P,gl)},\partial P)$ is diffeomorphic to $(N,T)$ through a diffeomorphism sending the framing of $\partial P$ in that of $T$.
\end{defi}
\begin{teo}[Turaev \cite{Tu}]\label{teo:existence}
Each pair $(N,T)$ as above admits a shadow.
\end{teo}
\begin{prf}{1}{
Let us sketch the idea of the proof, for more details see \cite{Tu} or \cite{Coin}. It is not difficult to check that $S^2$ equipped with gleam $1$ is a shadow $(P_0,gl_0)$ of $(S^3,\emptyset)$ (in particular $\partial P=\emptyset$) and that there is a collapse of $M_{(P_0,gl_0)}$ on $(P_0,gl_0)$ which restricts to a projection $\pi:S^3\rightarrow P_0$. Let $L$ be a framed link in $S^3$ such that $N$ is obtained by surgery on $L$; up to a small isotopy  we can suppose that $T$ does not intersect $L$ in $N$ and define $T'$ as the graph in $S^3-L$ corresponding to $T$. In $M_{(P_0,gl_0)}$ let $P_1$ be the polyhedron obtained by considering the mapping cylinder of the projection $\pi:L\cup T'\rightarrow P_0$; up to small isotopies of $L\cup T$ we can suppose it to be a simple polyhedron properly embedded in $M_0$ and whose boundary components form $L\cup T'$. By Proposition \ref{prop:gleam} we can equip $P_1$ with a gleam $gl_1$ so that $(P_1,gl_1)$ is a shadow of $(S^3,L\cup T')$. To conclude it is sufficient to glue to $P_1$ along each component of $L$ a disc equipped with zero gleam: indeed, at the level of the $4$-thickenings this corresponds to adding the $2$-handles which on the boundary produce $N$ out of $(S^3,L)$.
}\end{prf}
\subsection{Explict constructions}\label{sub:explconstr}
The proof of Theorem \ref{teo:existence} does not give a ``practical recipe" to construct shadows of $(N,T)$ because to find the gleams of a shadow of a trivalent graph in $S^3$ we used Proposition \ref{prop:gleam}. Hence we now describe a way of explicitly constructing a shadow of a framed graph $T$ in $S^3\#_k S^1\times S^2$. More precisely we produce shadows of $T\cup B$ where $B$ is an additional link we choose to add in order to get ``very simple" results; in our later applications, we will show that $B$ does not affect the calculation of the quantum invariants of $T$.  
Let us present $S^3\#_k S^1\times S^2$ as a surgery over $k$ $0$-framed unknots $u_i,\ i=1,\ldots k$ in $S^3$ and $T$ as a graph in $S^3$ whose edges are equipped with half-integers measuring the difference between the framing of $T$ around them and the blackboard-framing. Note that half-integers are necessary in general (see Remark \ref{rem:framing}).

\begin{figure}
\psfrag{T}{$T$}
\psfrag{B}{$B$}
\psfrag{0}{$0$}
\psfrag{-12}{$-\frac{1}{2}$}
\psfrag{12}{$\frac{1}{2}$}
\psfrag{-32}{$-\frac{3}{2}$}
\psfrag{1}{$1$}
 \centerline{\includegraphics[width=4.0cm]{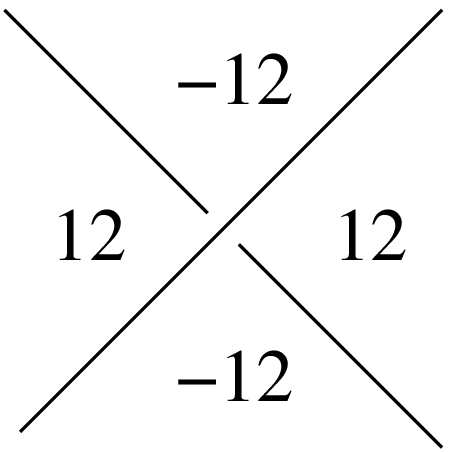} \includegraphics[width=5.0cm]{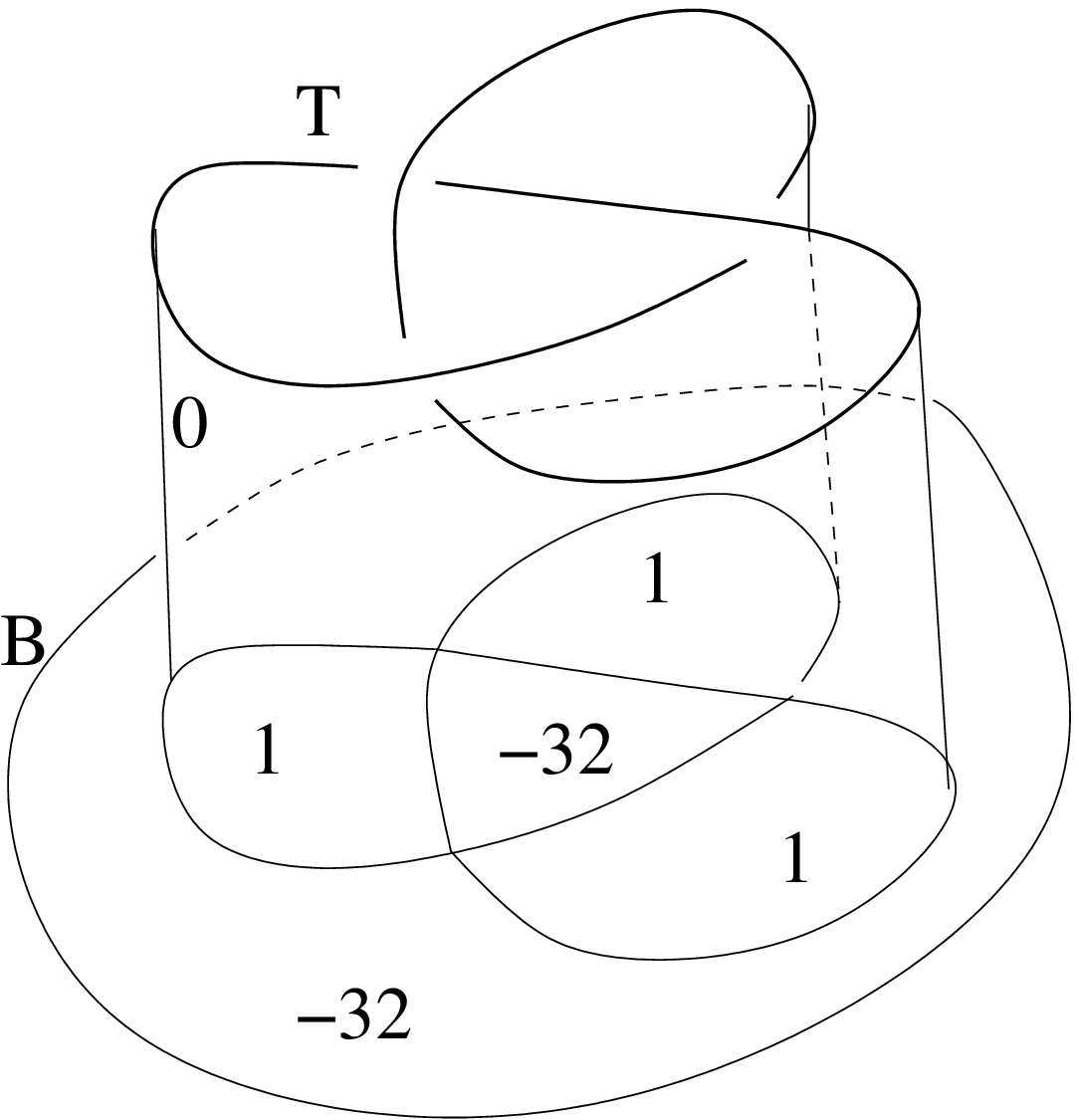}}
  \caption{In the left part of the figure, we show how to assign the local contributions to the gleams of the regions around a crossing. In the right, we exemplify the case of the trefoil. }
  \label{fig:pag408}
\end{figure}

Up to isotopy, we can suppose that a diagram $D$ of the projection in $\mr^2$ of $T\cup u_i,\ i=1,\ldots k$ is contained in a disc with $k-1$ holes $S_k$. The shadow we are searching for, is obtained by equipping with gleams the mapping cylinder $P$ of the projection of $T$ on $S_k$. To calculate the gleam of each region of $P$, let us note that each region either corresponds to a connected component of $S_k-D$ or to a cylinder over an edge of $T$. We equip the latter regions with the half integer written on the corresponding edge of $T$. To deal with the former regions, it is sufficient to sum up all the local half-integer contributions obtained by applying the rule of the left part of Figure \ref{fig:pag408}. The boundary of the resulting shadow $(P,gl)$ is composed of $T$ and of $k+1$ curves corresponding to $\partial S_k$; one can check that, applying Theorem \ref{teo:reconstruction}, the $4$-manifold $M_{(P,gl)}$ is a boundary connected sum of $B^4$ and $k$ copies of $S^1\times D^3$ and that $(\partial M_{(P,gl)}, \partial P)$ is diffeomorphic to $(S^3\#_k S^1\times S^2,T\cup B)$ where $B$ corresponds to $\partial S_k$. 
\begin{example}
Let $T$ be a left-handed trefoil in $S^3$ equipped with the blackboard framing in the diagram of Figure \ref{fig:pag408}. Following the above procedure as shown in the right part of the figure, one gets a polyhedron containing $3$ vertices (one per each crossing), $4$ disc regions equipped with gleams $1,1, 1,\ -\frac{3}{2}$ and two annular regions: the first, equipped with gleam $0$ and containing a free boundary component corresponding to $T$ and the second, equipped with gleam $-\frac{3}{2}$ containing a free boundary component corresponding to $B$. In that case $B\subset S^3$ is the unknot unlinked with $T$. 
\end{example}  
\section{The Jones invariant for framed links and trivalent graphs in $S^3\#_kS^2\times S^1$}\label{sec:jones}
In this section we start recalling how to calculate the Reshetikhin-Turaev invariants for a pair $(N,T)$ where $N$ is an oriented $3$-manifold and $T$ is a framed trivalent graph in $N$ using state-sums on shadows (see \cite{Tush} for a detailed account). Then we show how, when $N=S^3\#_kS^2\times S^1$, these state-sums can be used to produce invariants of pairs $(N,T)$ with values in the rational functions on $\mc$ which we will call $d$-colored Jones invariants. In the last subsection we prove some of the properties of these invariants and provide examples which will be crucial in the next section: in particular we recall the definition and the main properties of the universal hyperbolic links.
\subsection{State-sum quantum invariants}
Let $r,d$ be non-negative integers such that $r\geq 3$ and $1\leq d\leq r$, let $t$ be a complex valued variable and let 
for each $n\in \mathbb{N}$: $$[n]\doteq\frac{t^\frac{n}{2}-t^{-\frac{n}{2}}}{t^{\frac{1}{2}}-t^{-\frac{1}{2}}}$$ $$[n]!\doteq \prod_{1\leq i\leq n} [i],\ [0]!\doteq[1]!\doteq1$$ $$w_j\doteq (-1)^{2j}[2j+1]$$ 
We say that a triple $(i,j,k)$ of elements of $\frac{\mathbb{N}}{2}$ is \emph{admissible} if the following conditions are satisfied:
\begin{enumerate}
\item $|i+j|\leq k,\ |i+k|\leq j,\  |j+k|\leq i$;
\item $i+j+k \in\mathbb{N}$.
\end{enumerate}
For each admissible triple of elements of $\frac{\mathbb{N}}{2}$ let $$\Delta(i,j,k)\doteq \sqrt{\frac{[i+j-k]![i+k-j]![j+k-i]!}{[i+j+k+1]!}}$$
For any $6$-uple $(i,j,k,l,m,n)$ of elements of $\frac{\mathbb{N}}{2}$ such that the three-uples $(i,j,k),\ (i,m,n),\ (j,l,n)$ and $(k,l,m)$ are admissible, we define its \emph{$6j$-symbol} as follows:
$$ \setlength\arraycolsep{1pt}
\left( \begin{array}{ccc}
i & j &k\\
l  & m & n\\
\end{array}\right)\ns=\ns\sum_{z\in \mathbb{N}}\nns \frac{(\sqrt{-1})^{-2(i+j+k+l+m+n)}\Delta(i,j,k)\Delta(i,m,n)\Delta(j,l,n) \Delta(k,l,m) (-1)^z [z+1]!}{[z\nns -\nns i\nns -\nns j\nns -\nns k]![z\nns -\nns i\nns -\nns m\nns -\nns n]![z\nns-\nns j\nns -\nns l\nns -\nns n]![z\nns-\nns k\nns-\nns l\nns- \nns m]![ i\nns +\nns j\nns +\nns l\nns +\nns m\nns -\nns z]![i\nns +\nns k\nns +\nns l\nns +\nns m\nns -\nns z]![j\nns +\nns k\nns +\nns m\nns +\nns n\nns -\nns z]!}
$$
where the sum is taken over all $z$ such that all the arguments of the quantum factorials in the denominator of the r.h.s. are non-negative integers.
Let furthermore: $$W\doteq \sum_{0\leq i\leq \frac{r-2}{2}} (w_i)^2 \hspace{1cm }S\doteq W^{-1}\nns\nns\nns\sum_{0\leq i\leq \frac{r-2}{2}} (w_i)^2 t^{(-i(i+1))}(\sqrt{-1})^{4i}$$
\begin{defi}
An \emph{admissible coloring} on a shadow $(P,gl)$ is the assignment of a non-negative half-integer to each region of $P$ such that, for each edge of $Sing(P)$, the three-uple of colors associated to the three regions containing the edge is admissible. An $r$-admissible coloring is an admissible coloring such that the color of each region of $P$ does not exceed $\frac{r-2}{2}$ and such that for each edge of $Sing(P)$, the three-uple $(i,j,k)$ of colors on the regions containing the edge satisfy the additional condition $i+j+k\leq r-2$. A coloring is said to be \emph{relative to a coloring $d$ of $\partial P$} if each region containing an edge of $\partial P$ is colored with the same color as that edge. 
\end{defi}
Given a coloring on $P$, for each region $R$ let $w[R]\doteq w_j^{\chi(R)} t^{-j(j+1)gl(R)}(\sqrt{-1})^{4jgl(R)}$
where $j$ is the color of $R$, $gl(R)$ is its gleam and $\chi(R)$ is its Euler characteristic; similarly, to each vertex we associate its $6j$-symbol
where $(i,j,k,l,m,n)$ are the colors of the regions around the vertex and $(i,l)$ $(j,m)$ and $(k,n)$ are the pairs of colors corresponding to regions which, near the vertex, intersect only in the vertex itself. Finally, let $sign(P,gl)$ be the signature of the self intersection form of $H_2(M_{(P,gl)};\mz)$ and $nul(P,gl)$ be the dimension of the maximal real subspace of $H_2(M_{(P,gl)};\mr)$ contained in the annihilator of the form. 
The {\it $(r,d)$-state sum} of $(P,gl)$, denoted $|(P,gl)|^r_d$, is the following sum taken over all the $r$-admissible colorings of $P$ relative to the coloring on $\partial P$ whose value on each edge is $\frac{(d-1)}{2}$:
$$|(P,,gl)|^r_d\doteq \sum_{colorings}\prod_{regions} w[R]\prod_{vertices}\left( \begin{array}{ccc}
i & j &k\\
l  & m & n\\
\end{array}\right)$$  
\begin{rem}
For the sake of clarity, let us stress that the gleam of $P$ and the admissible colorings are two different and independent objects. The gleam of $P$ only comes into play inside the terms $w[R]$ associated to the regions of $P$.
\end{rem}
\begin{teo}[$(r,d)$-Reshetikhin-Turaev invariants \cite{Tush},\cite{Tu}]\label{teo:RT}
Let $N$ be a $3$-manifold, $T\subset N$ be a framed trivalent graph in $N$ and $(P,gl)$ be a shadow of $(N,T)$. The complex number $RT^r_d(N,T)\doteq W^{1-\chi(P) -nul(P,gl)}S^{-sign(P,gl)}|(P,gl)|^r_d(e^{\frac{2\pi\sqrt{-1}}{r}})$ is a homeomorphism invariant of $(N,T)$; in particular, it does not depend on the choice of a shadow of $(N,T)$.
\end{teo}
\begin{rem}
The factor $W^{1-\chi(P) -nul(P,gl)}S^{-sign(P,gl)}$ is to be considered as a ``normalization factor".
\end{rem}
The particular case when $T$ is a framed link in a $3$-manifold $N$ diffeomorphic to a $S^1$-bundle over a surface $F$ with Euler number $eul$, and $(P,gl)$ is obtained by projecting $T$ in the shadow $(F,eul)$ of $N$, was studied in detail by Turaev (\cite{Tush}). In that case, the above sum is simplified since the first two factors (the powers respectively of $W$ and of $S$) do not depend on $T$ and hence can be discarded. Moreover, when $N$ is $S^3$ and $(P.gl)$ is constructed as explained in Subsection \ref{sub:explconstr}, the following was proved inside the proof of Theorem 6.4 of the same paper:
\begin{teo}[\cite{Tush}]\label{teo:64}
Let $T$ be a framed link in $S^3$, $d$ and $r$ be integers as above, $Q=e^{\frac{2\pi\sqrt{-1}}{r}}$ and $J_d$ be the colored Jones polynomial of $T$ associated to the $d$-dimensional irreducible representation of $U_q(sl_2(\mc))$ normalized so that its value on the $0$-framed unknot is $[d]$. It holds: $ J_{d}(Q)=RT^{r}_d$.
\end{teo}

\subsection{The Jones invariant through shadows}
Note that, strictly speaking, Theorem \ref{teo:64} only states that the evaluation in some roots of unity of the colored Jones polynomials can be calculated through the state-sums $|(P,gl)|$, but it is well known that, in the case of links in $S^3$, the whole colored Jones polynomials can be calculated through shadow-based state-sums (\cite{KR}). In this subsection we extend this approach  to define a homeomorphism invariant $J_d$ of links (and trivalent graphs) in $S^3\#_k S^1\times S^2$ with values in the rational functions over $\mc$. 

Let $T$ be a framed link (or graph) in $N=S^3\#_k S^2\times S^1$ colored with the positive integer $d\geq 2$, and let $(P,gl)$ be a shadow of $(N,T)$ constructed as explained in Subsection \ref{sub:explconstr}; in particular $\partial P$ is composed by $T$ and some extra components forming a link $B$ and corresponding to the boundary of the surface $S_k$ used to construct $(P,gl)$. Observe that the set $Col_d$ of admissible colorings on $P$ satisfying the following conditions is finite:
\begin{enumerate}
\item the color of a region containing one component of $B$ is $0$;
\item the color of a region containing an edge of $T$ is $\frac{(d-1)}{2}$.
\end{enumerate}
Indeed by construction $P$ retracts on a graph and the following lemma applies:
\begin{lemma}\label{lem:retract}
A simple polyhedron has a finite number of admissible colorings relative to a coloring of the boundary if and only if it retracts on a graph.
\end{lemma}
\begin{prf}{1}{
Given half-integers $i,j$, the set of $k$ such that $(i,j,k)$ is admissible is finite; moreover if a polyhedron retracts on a graph there is at least one edge of $Sing(P)$ touched by two external regions. Then, the ``if" follows by an easy induction argument.
Viceversa if $P$ does not retract on a graph, then, after retracting it as long as possible, one gets a set of simple polyhedra connected along graphs. Then an infinite set of admissible colorings on $P$ can be constructed by extending the coloring whose value on each region of these polyhedra is $2k,\ k\geq 0$.
}\end{prf}

Let us then define the \emph{d-colored Jones invariant} of $(N,T)$ denoted $J_d(N,T)\in\mathbb{C}(t^\frac{1}{4})$ as follows: $$J_d(N,T)\doteq (-1)^{d-1} [d]^{b_1(N)-1}|(P,gl)|_d$$ where $|(P,gl)|_d$ is the above defined sum but taken over all the colorings of $Col_d$ and $b_1(N)$ is the first Betti number of $N$. Since it is a finite sum of rational functions (we are no longer fixing any value for $t$), its result is a rational function. 
The following holds: 
\begin{teo}\label{teo:jones}
The rational function $J_d(N,T)$ is a homeomorphism invariant of the pair $(N,T)$: if $(P',gl')$ is another shadow of $(N,T)$ which retracts on a graph, then the rational functions $|(P,gl)|^d$ and $|(P',gl')|^d$ are identical. The Jones invariants ``extend" the Reshetikhin-Turaev invariants of $(N,T)$: for each integer color $d\geq 2$ on $T$ there exists an integer $d_0(N,T)$ such that for each $r\geq d_0(N,T)$ it holds $RT_d^r(N,T)=(-1)^{d-1}W^{b_1(N)}[d]^{1-b_1(N)}J_d(N,T)(e^{\frac{2\pi\sqrt{-1}}{r}})$.
In the particular case when $N=S^3$ and $T$ is a framed link, then $J_d(N,T)$ is the Jones polynomial of $T$ normalized so that its value on the $0$-framed unknot is $1$. \end{teo}
\begin{prf}{1}{
Note that if $(P,gl)$ and $(P',gl')$ are shadows of $(N,T)$ which retract on graphs, then $\chi(P)=1-b_1(N)$, $nul(P,gl)=0=sign(P,gl)$ and the same equalities hold for $(P',gl')$.
The first statement is an immediate consequence of the second statement: indeed if $(P,gl)$ and $(P',gl')$ are two shadows of $(N,T)$, then the rational functions $|(P,gl)|_d$ and $|(P',gl')|_d$ coincide on all the roots of unity of degree greater than $d_0(N,T)$ and hence are identical.

To prove the second statement, let $d_0(N,T)$ be defined as follows:
$$d_0(N,T)\doteq 3\min_{(P,gl)}\max_{colorings}\max_{regions} color(R)$$
where $(P,gl)$ runs in the set of all shadows of $(N,T)$, the colorings run between all the admissible colorings of $(P,gl)$ extending the color $d$ on $T$ and $R$ runs between all the regions of $(P,gl)$. Let us now fix any shadow $(P,gl)$ of $(N,T)$: by the definition of $d_0(N,T)$ for each $r\geq d_0(N,T)$ the set of $(r,d)$-admissible colorings of $(P,gl)$ coincides with the set $Col_d$ of colorings extending the color $d$ on $T$, and so $|(P,gl)|^r_d(e^{\frac{2\pi\sqrt{-1}}{r}})=|(P,gl)|_d(e^{\frac{2\pi\sqrt{-1}}{r}})$. Hence, since $RT_d^r(N,T)=W^{1-\chi(P)-nul(P,gl)}S^{-sign(P,gl)}|(P,gl)|^r_d=W^{b_1(N)}S^0|(P,gl)|_d(e^{\frac{2\pi\sqrt{-1}}{r}})$ the thesis follows.

The last statement is immediate since we already know that RT-invariants are extended by the Jones polynomials (Theorem \ref{teo:64}) and, again, if two rational functions are equal on an infinite number of points, they are identical.
}
\end{prf}
\begin{rem}
Note that the above result states that $J_d(N,T)$ can be calculated through \emph {any} shadow which retracts on a graph and not only through a shadow constructed as explained in Subsection \ref{sub:explconstr}. This will be used to find easy formulas for $J_d(N,T)$ in some particular cases.
\end{rem}

\subsection{Properties and examples}\label{sub:examples}
Let $T$ be a framed link in $N=S^3\#_k S^2\times S^1$ and $T'$ be $T$ equipped with a framing which differs from that of $T$ by $s$ twists. The following holds:
\begin{lemma}\label{lem:framings}
For each $d\in \mathbb{N}, \ d\geq 2$, it holds $J_d(N,T')=(-1)^{s(d-1)}t^{-s\frac{d^2-1}{4}}J_d(N,T)$. 
\end{lemma}
\begin{prf}{1}{
It is sufficient to note that a shadow of $(N,T')$ can be obtained from a shadow $(P,gl)$ of $(N,T)$ by adding $s$ to the gleam of the region $R$ containing $T$. Since $R$ is by construction colored with $\frac{(d-1)}{2}$, this changes each summand of $|(P,gl)|_d$ by a constant factor equal to $(\sqrt{-1})^{4s\frac{d-1}{2}}t^{-s\frac{(d-1)}{2}\frac{(d+1)}{2}}$. 
}\end{prf}
The above lemma shows that $J_d(N,T)$ can be used to detect knots which have different but isotopic framings, as for instance, knots intersecting geometrically once an embedded sphere. Indeed the following holds:
\begin{cor}
If $k\subset N$ is a knot admitting two different framings which are isotopic, then $J_d(N,k)=0,\ \forall d$.  
\end{cor} 
\begin{example}
Let $(N,T)=(S^3,o)$ where $o$ is the $0$-framed unknot in $S^3$. Then a shadow of $(S^3,o)$ is given by a disc with gleam $0$. If $o$ is colored with the color $d$, then $|(P,gl)|_d=(-1)^{(d-1)}[d]$, and $J_d(S^3,o)=1$.
\end{example}
\begin{example}
Let $T$ be a left-handed trefoil in $S^3$ as in the right part of Figure \ref{fig:pag408} equipped with the blackboard framing and let $P$ be the polyhedron obtained by gluing a disc to the core of a Moebius strip. It is not difficult to check that, equipping the disc with gleam $-\frac{3}{2}$ and the annular region containing $\partial P$ with gleam $3$, $(P,gl)$ is a shadow of $(S^3,T)$: indeed $(P,gl)$ is obtained by collapsing the region which in Figure \ref{fig:pag408} contains the ``superfluous" boundary component $B$. Then it holds: $J_d(S^3,T)=(-1)^{(d-1)}[d]^{-1}t^{-3\frac{d^2-1}{4}}(\sqrt{-1})^{6(d-1)}\sum_{k\leq d-1}[2k+1]t^{\frac{3}{2}k(k+1)}(\sqrt{-1})^{-6k}$ with $k\in \mn$.
In particular, if $d=2$ one gets $J_2(S^3,T)=-[d]^{-1}t^{-\frac{9}{4}}(-1)(1-[3]t^{3})=t^{-\frac{9}{4}}\frac{1-t^{4}-t^{3}-t^{2}}{t^\frac{1}{2}+t^{-\frac{1}{2}}}=-t^\frac{9}{4}\frac{-1+t+t^3}{t^4}=-t^\frac{9}{4}J(T)$ where $J(T)$ is the ordinary Jones polynomial of the trefoil and the factor $-t^\frac{9}{4}$ comes from the fact that the Seifert framing of $T$ is twisted $-3$ times w.r.t. the framing encoded by the chosen shadow.
\end{example}

\begin{figure}
\psfrag{0}{$0$}
\psfrag{-12}{$-\frac{1}{2}$}
\psfrag{12}{$\frac{1}{2}$}
\psfrag{-32}{$-\frac{3}{2}$}
\psfrag{1}{$1$}
 \centerline{\includegraphics[width=10.0cm]{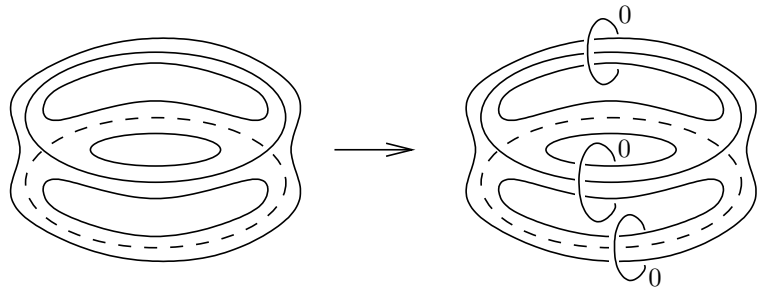}}
  \caption{In the left part of the figure, we show the regular neighborhood of the singular set of a simple polyhedron $P$: in that case $P$ has two vertices and $4$ edges. In the right, we exhibit a kirby diagram for $(N_P,L_P)$: in that case $L_P$ has $6$ connected components and the three $0$-surgered unknots represent the non-disconnecting spheres of $N_P=\#_3 S^2\times S^1$.  }
  \label{fig:exempleuniversal}
\end{figure}

\begin{example}
Let $T$ be the framed graph formed by the edges of a tetrahedron in $S^3$. A shadow of $T$ is given by the rightmost polyhedron in Figure \ref{fig:singularityinspine}, all whose regions are equipped with the gleam $0$. Equip all the edges of $T$ with the color $d$; if $d$ is even, then the set $Col_d$ is empty and so $J_d(S^3,T)=0$. Otherwise, it holds $$J_d(S^3,T)=[d]^{-1}[d]^6 \left( \begin{array}{ccc}
\frac{d-1}{2} & \frac{d-1}{2} &\frac{d-1}{2}\\
\frac{d-1}{2}  & \frac{d-1}{2} & \frac{d-1}{2}\\
\end{array}\right)$$
\end{example}
The above example shows that, in contrast with the case when $T$ is a link, $J_d(S^3,T)$ is not always a Laurent polynomial if $T$ is a graph. The same is true when $T$ is contained in $\#_k S^2\times S^1$, even if $T$ is a link. To show this, let us introduce a family of links contained in $\#_k S^2\times S^1$ called \emph{universal hyperbolic links}, which has been first studied by Dylan Thurston and the author (\cite{CT2}). Let $P$ be a simple polyhedron whose singular set is connected and containing $c(P)\geq 1$ vertices, and let $P'$ be the regular neighborhood of $Sing(P)$ in $P$, equipped with zero gleam. By Theorem \ref{teo:reconstruction}, $(P',0)$ can be thickened to a $4$-manifold $M_{(P',0)}$ which is diffeomorphic to a regular neighborhood of a graph in $\mr^4$ and whose boundary is $\#_{c+1} S^2\times S^1$; moreover $\partial P'\subset \partial M_{(P',0)}$ is a framed link. 
Let $(N_P,L_P)$ be the pair $(\partial M_{(P',0)},\partial P')$:
\begin{defi}[Universal Hyperbolic Link]\label{defi:unilink}
A link $L$ in $N=\#S^2\times S^1$ is said to be a \emph{universal hyperbolic link} if $(N,L)$ is diffeomorphic to $(N_P,L_P)$ for some simple polyhedron $P$ with $c(P)\geq 1$. 
\end{defi}
\begin{rem}
Different simple polyhedra could give the same pair: indeed only the combinatorics of a neighborhood of $Sing(P)$ in $P$ is relevant to the construction of $(N_P,L_P)$.
\end{rem}
\begin{example}\label{example:unilink} Let $(N_P,L_P)$ be a universal hyperbolic link and let $(P',0)$ the associated shadow whose complexity is $c$. If $d$ is an even integer then the coloring of $\partial P'$ given by $\frac{d-1}{2}$ does not extend to any coloring of $P'$, and hence $J_d(\partial M_{(P',0)},\partial P')=0$. If $d$ is odd, since each region of $P'$ has zero Euler characteristic and zero gleam, it holds:
$$J_d(\partial M_{(P',0)},\partial P')=[d]^{c} \left( \begin{array}{ccc}
\frac{d-1}{2} & \frac{d-1}{2} &\frac{d-1}{2}\\
\frac{d-1}{2}  & \frac{d-1}{2} & \frac{d-1}{2}\\
\end{array}\right)^c$$
\end{example}
The following theorem summarizes some of the most important properties of the universal hyperbolic links:
\begin{teo}[F.Costantino-D.P.Thurston \cite{CT2}]
Let $(N_P,L_P)$ be a universal hyperbolic link associated to a polyhedron $P$ whose complexity is $c(P)$. The complement of $L_P$ in $N_P$ can be equipped with a complete, hyperbolic metric with volume $2c(P)Vol_{Oct}$, where $Vol_{Oct}$ is the volume of an ideal regular octahedron in $\mathbb{H}^3$. 

Each orientable $3$-manifold $Q$ with boundary composed by tori or empty, can be obtained as an integer Dehn filling over a link $L_P$ for a suitable $P$. Moreover, if $Q$ is hyperbolic, $L_P$ can be chosen so that the following holds:
$$k\sqrt{c(P)} \leq Vol(Q)\leq 2Vol_{Oct}c(P)$$ 
where $k$ is a positive constant which does not depend on $Q$.
\end{teo}
We stress here that the procedure which associates a link to a simple polyhedron is constructive. An example is given in Figure \ref{fig:exempleuniversal}.

\section{The Volume Conjecture for links in $\#_k S^2\times S^1$}
In this section we first discuss how to extend the Volume Conjecture to include the case of links in $\#_k S^2\times S^1$. Then, we prove that the conjecture is true for all the universal hyperbolic links (see Definition \ref{defi:unilink}). 
Let $L$ be a framed link in $S^3$ and let $J_d(L)(t)$ be the d-colored Jones polynomial of $L$ (viewed as a variable in $t$) normalized so that its value on the unknot is $1$. The Volume Conjecture states the following:
$$\lim_{d\rightarrow \infty} log(|J_d(L)(e^{\frac{2\pi\sqrt{-1}}{d}})|)=\frac{Vol(L)}{2\pi}$$
where $Vol(L)$ is the hyperbolic volume of the complement of $L$ in $S^3$ (or its Gromov norm if $L$ is not hyperbolic).
In order to extend the VC to the setting of framed links in $\#_k S^2\times S^1$, three remarks are in order:
\begin{enumerate}
\item The invariant $J_d$ of framed links in $\#_k S^2\times S^1$ depends on a framing but, by Lemma \ref{lem:framings}, the modulus of its evaluation on a root of $1$ does not depend on it.
\item For some link $T$ there could exist colorings admitting no extension to a shadow of $T$ (see Example \ref{example:unilink}) and hence for which $J_d(N,T)=0$.
\item For some link $T$ the rational function $J_d(N,T)$ could have a pole at the $d$-th root of unity. 
\end{enumerate}
Hence, we choose to formulate the following extension of the VC using the Jones invariants defined in Section \ref{sec:jones} as follows:
\begin{conj}[Extended Volume Conjecture]
Let $T$ be a framed link in $N=\#_k S^2\times S^1$ and let $J_d(N,T)(t)$ be its d-colored Jones invariant. It holds:
$$\lim_{d\rightarrow \infty} \frac{1}{d}log(|J_d(N,T)(e^{\frac{2\pi\sqrt{-1}}{d}})|)=\frac{Vol(T)}{2\pi}$$
where $Vol(T)$ is the hyperbolic volume of $N-T$ and the limit is taken over all the $d\in \mathbb{N}$ such that $J_d(N,T)(t)\neq 0$ and has no pole in $e^{\frac{2\pi\sqrt{-1}}{d}}$.
\end{conj}
\begin{rem}
We use the notation ``Extended Volume Conjecture" to distinguish it from the ``Generalized Volume Conjecture" (see \cite{GL}), which deals with the asymptotic behavior of the evaluation of the Jones invariants in $e^{\frac{2\pi \sqrt{-1}\alpha}{d}}$, with $\alpha\in (0,1)$. 
\end{rem}

\begin{teo}\label{teo:vcvera}
The Extended Volume Conjecture is true for all universal hyperbolic links.
\end{teo}
Before proving the above result, let us recall some definitions and facts.
For each $x\in \mathbb{R}$ let us define the Lobatchevskji function $\Lambda(x)=-\int_0^x log(|2sin(s)|)ds$; $\Lambda(x)$ is smooth out of $\{ \pi k,\ k \in \mz\}$ and $\pi$-periodic. The Lobatchevskji function is a crucial element to calculate the volume of an ideal hyperbolic simplexes: for instance, the volume of the ideal regular octahedron is $Vol_{Oct}=8\Lambda(\frac{\pi}{4})$.
\begin{figure}
\psfrag{yl}{$\Lambda(x)$}
\psfrag{xl}{$x$}
\psfrag{0.5}{$0.5$}
\psfrag{1}{$1$}
\psfrag{1.5}{$1.5$}
\psfrag{2}{$2$}
\psfrag{2.5}{$2.5$}
\psfrag{3}{$3$}
\psfrag{-0.4}{$-0.4$}
\psfrag{-0.2}{$-0.2$}
\psfrag{0.4}{$0.4$}
\psfrag{0.2}{$0.2$}
 \centerline{\includegraphics[width=10.0cm]{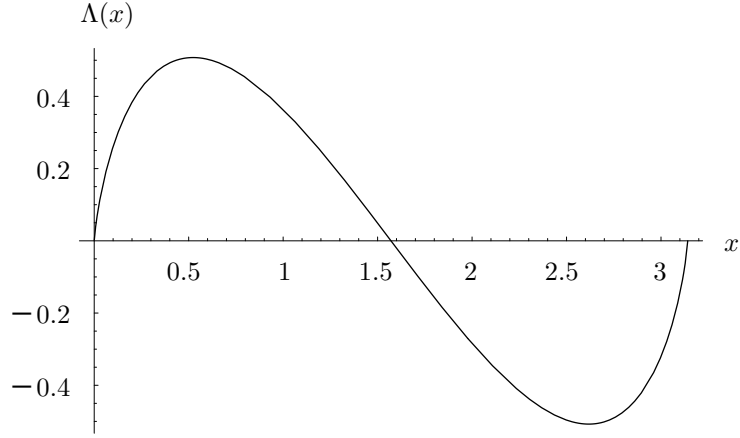}}
  \caption{A plot of $\Lambda (x), \ x\in [0,\pi]$.  }
  \label{fig:lambda}
\end{figure}

Given a rational function $r(t)$, let us call $ev_n(r)$ the evaluation of $r$ in $e^{\frac{2\pi\sqrt{-1}}{n}}$. Moreover, for each integers $l\geq s$
let us define the quantum binomials as:
$$\left[ \begin{array}{c}
l \\
s\\
\end{array}
\right] \doteq \frac{[l]!}{[s]![l-s]!}$$
\begin{lemma}\label{lem:legarouf}
There exists a constant $C$ such that for each pair of integers $(k,j)$ with $k\geq 2$ and $1\leq j\leq k-1$, letting $d=2k+1$ it holds:
$$e^{-\frac{d}{\pi}\Lambda(\pi\frac{j}{d})-Clog(d)}\leq (2sin(\frac{\pi}{d}))^jev_d([j]!)\leq e^{-\frac{d}{\pi}\Lambda(\pi\frac{j}{d})+Clog(d)}\\$$
$$e^{-\frac{d}{\pi}(\Lambda(\pi\frac{k}{d})-\Lambda(\pi\frac{j}{d})-\Lambda(\pi\frac{k-j}{d}))-3Clog(d)}\leq ev_d(\left[ \begin{array}{c}
k \\
j\\
\end{array}
\right])\leq e^{-\frac{d}{\pi}(\Lambda(\pi\frac{k}{d})-\Lambda(\pi\frac{j}{d})-\Lambda(\pi\frac{k-j}{d}))+3Clog(d)}$$ 
\end{lemma}
\begin{prf}{1}{
The existence of $C$ satisfying the former inequalities was proven by Garoufalidis and Le (\cite{GL}, Lemma 4.1): note indeed that our definition of $[j]$ differs from that in \cite{GL} by a factor $2sin(\frac{\pi}{d})$. We stress here that $C$ does not depend neither on $k$ nor on $i$. The latter inequalities are a straightforward consequence of the definition of the quantum binomial and of the first inequalities.
}
\end{prf}
\begin{lemma}\label{lem:mioanalitico}
Let $L(\alpha,\beta)\doteq \Lambda(\alpha)+\Lambda(\frac{\pi}{2}-\alpha+\beta)-\Lambda(\frac{\pi}{2}+\beta)$. For each interval $[a,b]\subseteq[0,\frac{\pi}{2}]$, let  $m_{[a,b]}$ and $M_{[a,b]}$ be respectively the min and the max in $[a,b]$ of $\Lambda(x)+\Lambda(\frac{\pi}{2}-x)$. For each $\epsilon > 0$ it holds:
\begin{enumerate}
\item There exists $\delta(\epsilon)$ such that $m_{[a,b]}-\epsilon<L(\alpha,\beta)<M_{[a,b]}+\epsilon,\ \forall (\alpha,\beta)\in [a,b]\times[-\delta(\epsilon),\delta(\epsilon)]$.
\item There exists $\delta(\epsilon)$ such that $L(\alpha,\beta)<2\Lambda(\frac{\pi}{4})+\epsilon,\ \forall (\alpha,\beta)\  s.t. \ |\beta|<\delta$.
\end{enumerate}
\end{lemma}
\begin{prf}{1}{
The first statement is implied by the fact that $\Lambda(x)$ and $L(\alpha,\delta)$ are continuous and that $L(\alpha,0)=\Lambda(x)+\Lambda(\frac{\pi}{2}-x)$.
The second statement is an instance of the first one when $[a,b]=[0,\frac{\pi}{2}]$. In that case $M_{[a,b]}=\max_{[0,\frac{\pi}{2}]} \Lambda(x)+\Lambda(\frac{\pi}{2}-x)=2\Lambda(\frac{\pi}{4})$ since $\Lambda'(x)+\Lambda'(\frac{\pi}{2}-x)=0$ is satisfied only when $x=\frac{\pi}{4}$.
}\end{prf}
\begin{prf}{3}{
Let $(N_P,L_P)$ be a universal hyperbolic link associated to a polyhedron $P$ containing $c$ vertices. For each $k\in \mathbb{N}$, let $d=2k+1$.
By Example \ref{example:unilink} we need to prove:
$$\lim_{k\rightarrow \infty} \frac{1}{d} log(|ev_{d}([d]^c \left( \begin{array}{ccc}
k & k & k\\
k  & k & k\\
\end{array}\right)^c)|)=2c\frac{Vol_{Oct}}{2\pi}=16c\frac{\Lambda(\frac{\pi}{4})}{2\pi}$$
It is sufficient to prove the above limit for $c=1$. Now, to simplify $ev_{d}(J_{d}(N,T))$, let us note that, for each $i\leq 2k$, the following equalities hold:
$$ev_{d}([i])=ev_{d}([2k+1-i])=-ev_{d}([2k+1+i]),\ ev_{d}([i]!)=ev_d(\frac{[2k]!}{[2k-i]!})$$
Then we apply them to simplify $J_d(N,L)$:
$$ev_d(J_d(N,L))=[d](\frac{[k]![k]![k]!}{[3k+1]!})^2\sum_{j=0}^{j=k} (-1)^{(3k+j)}\frac{[3k+1+j]!}{[j]!^4[k-j]^3}=\\\frac{1}{[2k+1]}\sum_{j=0}^{j=k}(-1)^{4k+2j}\frac{[2k+1][2k]!^2}{[j]!^4[k-j]!^4}=$$
$$=\sum_{j=0}^{j=k} (\frac{[k]!}{[j]![k-j]!})^4=\sum_{j=0}^{j=k}\left[ \begin{array}{c} k\\ j\\ \end{array}\right]^4$$
Where all the above equalities hold when evaluated at $e^{\frac{2\pi\sqrt{-1}}{d}}$. Note that the summands in the r.h.s are positive real numbers.

We first claim that $\limsup_{d\to \infty}\frac{1}{d}log(|ev_d(J_d(N,L))|)\leq \frac{8}{\pi}\Lambda(\frac{\pi}{4})$. By Lemma \ref{lem:legarouf} and point $2$ of Lemma \ref{lem:mioanalitico}, there exists a constant $C$ such that: 
$$ev_d(J_d(N,L))= \left[ \begin{array}{c} k\\ 0\\ \end{array}\right]^4+\left[ \begin{array}{c} k\\ k\\ \end{array}\right]^4+\sum_{j=1}^{j=k-1}\left[ \begin{array}{c} k\\ j\\ \end{array}\right]^4\leq 2+ (k-1)(e^{-\frac{d}{\pi}(2\Lambda(\frac{\pi}{4})+\epsilon)+3Clog(d)})^4$$
Where $\epsilon>0$ is any positive constant and $k$ is big enough so that $|\frac{k}{2k+1}-\frac{1}{2}|<\delta(\epsilon)$. Then, taking the logarithm and dividing by $d$, one obtains that $\limsup_{d\to \infty} \frac{1}{d}log(|ev_d(J_d(N,L))|)$ is not greater than $\Lambda(\frac{\pi}{4})\frac{8}{\pi}+\frac{4}{\pi}\epsilon$. But since $\epsilon$ can be chosen arbitrarily small the claim follows.

We now prove that the limit indeed exists and that it equals $\Lambda(\frac{\pi}{4})\frac{8}{\pi}$: to do this, we show that $\liminf_{d\to \infty} \frac{1}{d}log(|ev_d(J_d(N,L))|)\geq\Lambda(\frac{\pi}{4})\frac{8}{\pi}$. 

Let $\alpha \in (0,1)$ and let us bound $ev_d(J_d(N,L))$ from below as follows:
$$ev_d(J_d(N,L))\geq \sum_{j\geq\frac{\alpha k}{2}}^{j\leq \frac{k}{2\alpha}}\left[ \begin{array}{c} k\\ j\\ \end{array}\right]^4$$
Then, we now apply Lemma \ref{lem:legarouf} and 1) of Lemma \ref{lem:mioanalitico} with $[a,b]=[\frac{\alpha\pi}{4},\frac{\pi}{4\alpha}]$ and $m(\alpha)=m_{[\frac{\alpha k}{2},\frac{k}{2\alpha}]}$. It holds:
$$\sum_{j\geq\frac{\alpha k}{2}}^{j\leq \frac{k}{\alpha 2}}\left[ \begin{array}{c} k\\ j\\ \end{array}\right]^4\geq \sum_{j\geq \frac{\alpha k}{2}}^{j\leq \frac{k}{2\alpha}}(e^{\frac{d}{\pi}(m(\alpha)-\epsilon(\frac{k}{d}-\frac{1}{2}))-3Clog(d)})^4\geq \frac{k}{2}(\frac{1}{\alpha}-\alpha)e^{\frac{4d}{\pi}(m(\alpha)-\epsilon(\frac{1}{2}-\frac{k}{d}))-12Clog(d)}$$
where $\lim_{d\to \infty}\epsilon(\frac{1}{2}-\frac{k}{d})=0$.
Taking the logarithm and dividing by $d$ we get:
$$\frac{1}{d}log(|ev_d(J_d(N,L))|)\geq \frac{4}{\pi}(m(\alpha)-\epsilon(\frac{k}{d}-\frac{1}{2}))-12C\frac{log(d)}{d}$$ 
Then, letting $d\to \infty$ we get:
$$\liminf_{d\to \infty}\frac{1}{d}log(|ev_d(J_d(N,L))|)\geq \frac{4}{\pi}m(\alpha),\ \forall \alpha\in (0,1)$$
But then, since $\lim_{\alpha \to 1} m(\alpha)=2\Lambda(\frac{\pi}{4})$ the thesis follows. 
}\end{prf}

\begin{figure}
\psfrag{xdati}{$d$}
\psfrag{ydati}{$\frac{1}{d} log(|ev_d(J_d)|)$}
\psfrag{100}{$100$}
\psfrag{200}{$200$}
\psfrag{300}{$300$}
\psfrag{400}{$400$}
\psfrag{500}{$500$}
\psfrag{1.04}{$1.04$}
\psfrag{1.06}{$1.06$}
\psfrag{1.08}{$1.08$}
\psfrag{1.12}{$1.12$}
\psfrag{1.14}{$1.14$}
 \centerline{\includegraphics[width=10.0cm]{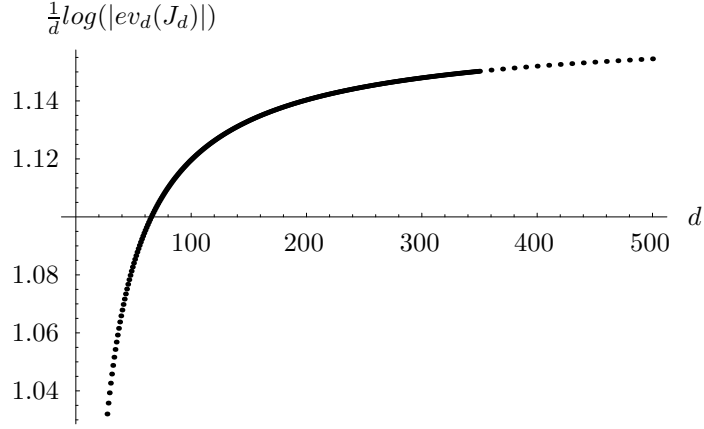}}
  \caption{A plot of $\frac{1}{d} Log(|J_d(N_P,L_P)|)$ for a universal hyperbolic link with complexity $1$: the limit is $\frac{8}{\pi}\Lambda(\frac{\pi}{4})=1.16624...$. }
  \label{fig:graficojones}
\end{figure}

We conclude by noting that it can be proved through algebraic arguments that the succession $\frac{1}{d} log(|J_d(N_P,L_P)|)$ is monotone increasing. Instead of proving it we limit ourselves to exhibit a computer-plot of the first $500$ elements of the succession (see Figure \ref{fig:graficojones}) .

\noindent

\end{document}